\documentstyle[amssymb,12pt]{amsart}

\newcommand{\labell}[1] {\label{#1}}

\numberwithin{equation}{section}

\newtheorem{thm}{Theorem}
\newtheorem{lemma}[thm]{Lemma}

\newtheorem{cor}[thm]{Corollary}

\theoremstyle{definition}

\newtheorem{rmk}[thm]{Remark}

\expandafter\chardef\csname pre amssym.def at\endcsname=\the\catcode`\@ 
\catcode`\@=11 
\def\undefine#1{\let#1\undefined} 
\def\newsymbol#1#2#3#4#5{\let\next@\relax 
 \ifnum#2=\@ne\let\next@\msafam@\else 
 \ifnum#2=\tw@\let\next@\msbfam@\fi\fi 
 \mathchardef#1="#3\next@#4#5}
\def\mathhexbox@#1#2#3{\relax 
 \ifmmode\mathpalette{}{\m@th\mathchar"#1#2#3}%
 \else\leavevmode\hbox{$\m@th\mathchar"#1#2#3$}\fi} 
\def\hexnumber@#1{\ifcase#1 0\or 1\or 2\or 3\or 4\or 5\or 6\or 7\or 8\or 
 9\or A\or B\or C\or D\or E\or F\fi} 
 
\font\teneufm=eufm10 
\font\seveneufm=eufm7
\font\fiveeufm=eufm5 
\newfam\eufmfam
\textfont\eufmfam=\teneufm 
\scriptfont\eufmfam=\seveneufm 
\scriptscriptfont\eufmfam=\fiveeufm 
 
\catcode`\@=\csname pre amssym.def at\endcsname 

\newcommand{\supp}{{\mathit supp}}

\newcommand{\const}{{\mathit const}}

\newcommand{\hook}{\hookrightarrow}
\newcommand{\tit}{\tilde{t}}
\newcommand{\ty}{\tilde{y}}

\def	\reals	{{\Bbb R}}

\def	\rationals	{{\Bbb Q}}
\def	\TT	{{\Bbb T}}

\def	\p	{\partial}

\def\qed{\smallskip\noindent
$\Box$ 
\smallskip}

\begin{document}


\setlength{\smallskipamount}{6pt}
\setlength{\medskipamount}{10pt}
\setlength{\bigskipamount}{16pt}





\title[The Hamiltonian Seifert conjecture in $\reals^6$]{A smooth 
counterexample to the Hamiltonian Seifert conjecture in $\reals^6$}

\author[Viktor Ginzburg]{Viktor L. Ginzburg}

\address{Department of Mathematics, UC Santa Cruz, 
Santa Cruz, CA 95064}
\email{ginzburg@@cats.ucsc.edu}

\date{May, 1997; Preprint: dg-ga/9703006, revised version}

\thanks{The work is partially supported by the NSF and by the faculty
research funds of the University of California, Santa Cruz.}

\bigskip

\begin{abstract}
A smooth counterexample to the Hamiltonian Seifert conjecture in $\reals^6$
is found. Namely, we construct a smooth proper function 
$H\colon\reals^{2n}\to\reals$, $2n\geq 6$, such that the level set $\{H=1\}$ 
is non-singular and has no periodic orbits for the Hamiltonian flow of $H$
with respect to the standard symplectic structure. The function
$H$ can be taken to be $C^0$-close and isotopic to a positive-definite
quadratic form so that $\{H=1\}$ is isotopic to an ellipsoid. This is a 
refinement of previously known constructions giving such
functions only for $2n\geq 8$. The proof is based on a 
new version of a symplectic embedding theorem applied to the horocycle flow.
\end{abstract}

\maketitle

\section{Introduction} 

Hamiltonian systems tend to have periodic orbits whenever the energy
level is non-singular and compact. For example, Hofer, Zehnder, and Struwe 
(see \cite{ho-ze:per-sol}, \cite{str}, and \cite{ho-ze:book})
have shown that
periodic orbits must exist on almost all energy levels of a
proper smooth function on the standard symplectic $\reals^{2n}$. 
This implies Viterbo's theorem, \cite{vi:thm}, proving 
Weinstein's conjecture, \cite{we:conj}:
a compact smooth hypersurface of contact type in the standard $\reals^{2n}$
carries at least one closed characteristic. (The reader interested in
a detailed discussion should consult \cite{ho-ze:book}.)

The first examples of smooth compact hypersurfaces
in $\reals^{2n}$ without closed characteristics were found only recently. 
As a consequence, one obtains a smooth proper function on $\reals^{2n}$
which has at least one regular level set without periodic orbits of its
Hamiltonian flow. These examples, constructed in \cite{gi:seifert}
and, independently, by M. Herman \cite{herman-fax} 
required the ambient space to be at least eight-dimensional, i.e., 
$2n\geq 8$, in the smooth case. A $C^{3-\epsilon}$-smooth hypersurface 
in $\reals^6$ was found by M. Herman \cite{herman-fax}. In the present 
paper we extend the results of \cite{gi:seifert} to $2n=6$ by 
constructing such examples in dimension six (Section \ref{sec:main}).

To carry out this extension, we show that,
as a consequence of a general symplectic embedding theorem proved below,
the horocycle flow can
be ``symplectically'' embedded into $\reals^5$. (See Section \ref{sec:proofs}.)
This is the only part of the argument of \cite{gi:seifert} missing
in the six-dimensional case. The paper is concluded by a list of all
Hamiltonian flows known to the author to have no periodic orbits 
(Section \ref{sec:list}).

Although formally independent, this paper can be considered as a 
follow-up of \cite{gi:seifert}. In particular, we refer to some of the
proofs, not only the results, given therein.
When possible and appropriate the original notation is kept.

\bigskip
\noindent{\bf Acknowledgments.} The author is deeply grateful to
Yasha Eliashberg and Richard Montgomery for useful discussions.

\section{Main Theorems} 
\labell{sec:main}

Let $i\colon Q\hookrightarrow V$ be an embedded smooth compact hypersurface
without boundary in a $2n$-dimensional symplectic manifold $(V,\sigma)$. 
Recall that a characteristic of a two-form $\eta$ of rank
$(2n-2)$ on a $(2n-1)$-dimensional manifold is an integral curve
of the field of directions formed by the null-spaces $\ker\eta$.
Our main result is the following 

\begin{thm}
\labell{thm:main1}
Assume that $2n\geq 6$ and $i^*\sigma$ has only a finite number of 
closed characteristics. Then there exists a $C^\infty$-smooth embedding 
$i'\colon Q\to V$, which is $C^0$-close and isotopic to $i$, such that 
${i'}^*\sigma$ has no closed characteristics.
\end{thm}

\begin{rmk}
For $2n\geq 8$, Theorem \ref{thm:main1} as well as other results stated in
this section are known. See \cite{gi:seifert} and
\cite{herman-fax}. This paper improves the dimensional 
constraint by two. As mentioned above, for $2n=6$, an example of such
a  $C^{3-\epsilon}$-smooth embedding $i'$ was found by
M. Herman \cite{herman-fax}.

As is the case for $2n\geq 8$, the embedding $i'$ can be chosen to
coincide with $i$ outside a finite number of small balls each ``centered''
on a closed characteristic of $i^*\sigma$.
\end{rmk}

An irrational ellipsoid $Q$ in the standard symplectic vector space 
$\reals^{2n}=V$ corresponds to a collection of $n$ harmonic
oscillators (a quadratic Hamiltonian) whose frequencies are linearly 
independent over $\rationals$. Thus $Q$ has exactly $n$ periodic orbits. 
Applying Theorem \ref{thm:main1}, we obtain

\begin{cor}    
\labell{cor:sphere}
For $2n\geq 6$, there exists a $C^\infty$-smooth embedding 
$S^{2n-1}\to \reals^{2n}$
such that the pull back of $\sigma_{2n}$ to $S^{2n-1}$ has no
closed characteristics.
\end{cor}
\begin{cor}
\labell{cor:function}
For $2n\geq 6$, there exists a $C^\infty$-smooth function 
$h\colon\reals^{2n}\to \reals$, $C^0$-close and isotopic
(with a compact support) to a positive definite quadratic form,
such that the Hamiltonian flow of $h$ has no closed trajectories on 
the level set $\{ h=1\}$.
\end{cor}

A similar result involving no symplectic embeddings, but only
two-forms on $Q$, is used in the proof of Theorem \ref{thm:main1}.
Let $\dim Q=2n-1$
and let $\eta$ be a closed maximally non-degenerate (i.e., of rank
$(2n-2)$) two-form on $Q$.

\begin{thm}
\labell{thm:main2}
Assume that $2n-1\geq 5$ and that $\eta$ has a finite number of 
closed characteristics. Then there exists closed a maximally non-degenerate
2-form $\eta'$ on $Q$ (homotopic to $\eta$) without closed characteristics.

\end{thm}
Recall that the forms $\eta$ and $\eta'$ are said to be {\em homotopic} 
if there exists a family $\eta_\tau$, $\tau\in [0,1]$, 
of closed maximally non-degenerate forms in the same cohomology class
connecting $\eta=\eta_0$ with $\eta'=\eta_1$.

\begin{rmk}
Theorem \ref{thm:main2} extends to the real analytic case: one can
make the form $\eta'$ real analytic, provided that $Q$ and $\eta$ are 
real analytic.
The argument is the same as that used in the construction of a real analytic
version of Wilson's, \cite{wilson}, or Kuperberg's, \cite{kuk}, flow.
(See \cite{ghys}.)
\end{rmk}

\section{Proofs} 
\labell{sec:proofs}

As is pointed out in \cite{gi:seifert} (Remark 3.5), to prove Theorems 
\ref{thm:main1} and \ref{thm:main2} it
is sufficient to find a ``symplectic'' embedding of the horocycle
flow into $\reals^5$. (We will elaborate on this later.) The existence
of such an embedding is a consequence of the following general construction.

Let $N$ and $W$ be manifolds of
equal dimensions. The manifold $N$ is assumed to be compact, perhaps with
boundary, while $W$ may be open but must be a manifold without 
boundary. Let $\sigma$ be a symplectic form on $W$. Abusing notation,
also denote by $\sigma$ the pull-back of $\sigma$ to $W\times\reals$ 
under the natural projection $W\times\reals\to W$. To avoid confusion,
we will sometimes indicate the domain of a form by a subscript, e.g., 
$\sigma_W$ or $\sigma_{W\times\reals}$.

\begin{thm}
\labell{thm:embed}
Let $\omega_t$, $t\in [0,1]$, be a family of symplectic forms
on $N$ in a fixed cohomology class: $[\omega_t]=\const$. Assume
that there is an embedding $j_0\colon N\to W\times\reals$ such
that $j_0^*\sigma=\omega_0$. Then there exists an embedding 
$j_1\colon N\to W\times\reals$, isotopic to $j_0$, with
$j_1^*\sigma=\omega_1$.
\end{thm}

\begin{rmk}
\labell{rmk:moser1}
Since $\sigma_W$ is symplectic, the composition of $j_1$ with
the projection to $W$ is necessarily an immersion.
When $\p N=\emptyset$, Theorem \ref{thm:embed} follows immediately
from Moser's theorem \cite{moser}. 
Here, however, we are more interested in
the case where $\p N\neq \emptyset$ and Moser's theorem does
not apply.
\end{rmk}

\begin{pf} Assume first that $N$ is a manifold with boundary.

Following Gromov, \cite{gr:book} (p. 336), let us approximate the
family $\omega_t$ by a family which piecewise linearly interpolates
a finite sequence of symplectic forms $\omega_{(k)}$, $k=0,\,1,\ldots, q$, 
such that 
$\omega_{(k)}=\omega_{(k-1)}+df_k\wedge dg _k$, where $f_k$ and $g_k$
are smooth functions on $N$. We assume, of course, that the new and
the original families have the same end-points: $\omega_0=\omega_{(0)}$ and
$\omega_1=\omega_{(q)}$. The existence of this approximation is easy
to verify.

Arguing inductively, we will find an embedding $i_k$ with
$i_k^*\sigma=\omega_{(k)}$ using the existence of $i_{k-1}$ with
$i_{k-1}^*\sigma=\omega_{(k-1)}$. The base of induction is given by
$i_0=j_0$. After $i_k$ is constructed, it will remain to set
$j_1=i_q$.

Thus assume that we have $i_{k-1}$. For the sake of simplicity, we
identify $N$ with $i_{k-1}(N)$ and denote $f_k$ and $g_k$ by
$f$ and $g$. The construction of $i_k$ will be carried out in two
steps: first, find a symplectic embedding 
$\tilde{i}\colon N\to W\times\reals^2$ 
and then push it back into $W\times\reals$.

{\em The first step.}  
Denote by $(t,y)$ the standard linear coordinates
on $\reals^2$. Let us think of $W\times\reals$ as given by the equation 
$y=0$ in $W\times\reals^2$ so that  the form $\sigma_{W\times\reals}$ is the 
restriction of $\sigma_W+dt\wedge dy$. However, to construct $\tilde{i}$, it
is more convenient to use a different system of coordinates on a
neighborhood of $N$. 

First, let us extend $N$ a little bit beyond its boundary so that
the resulting open manifold $N_1$ containing $N$ is still embedded 
into $W\times\reals$ transversely to the direction of the $t$-axis. 
As a consequence,
the restriction of $\sigma_W$ to $N_1$ is a symplectic form, which
we denote for the sake of simplicity by $\sigma$. Note that
the restriction of $\sigma$ to $N$ is just $\omega_{(k-1)}$, the
same as $\sigma_{W\times\reals}|_N$. In what follows we will need to take $N_1$ so that
the closure $\bar{N}_1$ is a smooth embedded manifold with boundary and
that $\bar{N}_1$ is still transversal to the $t$-axis. 
Since $W\times\reals$ is embedded into $W\times\reals^2$, so are
$N_1$ and $\bar{N}_1$.

Let $I=(-\epsilon,\epsilon)$ and $J=(-\delta,\delta)$ and
let $\tilde{t}$ and $\tilde{y}$ be the natural coordinates on
these intervals.

\begin{lemma}
For sufficiently small $\epsilon>0$ and $\delta>0$,
there exists an embedding
$\psi\colon N_1\times I\times J\to W\times\reals^2$ such that
$$
\psi^*(\sigma_W+dt\wedge dy) =\sigma+d\tilde{t}\wedge d\tilde{y}
$$
and, in addition, $\psi|_{N_1}$ is the original embedding 
$N_1\hook W\times \reals$ and 
$\psi(N_1\times I\times 0)\subset W\times\reals$.
\end{lemma}

The existence of $\psi$ is a consequence of the symplectic neighborhood
theorem applied to a small neighborhood of $\bar{N}_1$ in $W\times\reals$
considered as a submanifold in $W\times\reals^2$. For the sake of
completeness we give a direct proof.

\smallskip

\noindent{\em Proof} of the lemma.
To find $\psi$, let us first extend $\bar{N}_1$ beyond its boundary
in the same way as $N_1$ extends $N$. In particular, the closure of the
resulting manifold $N_2$ is assumed to be transversal to the $t$-axis
direction. Let, as before, $\sigma=\sigma_{W\times\reals}|_{N_2}$.

A neighborhood $U$ of $N_2$ in $W\times\reals^2$ is diffeomorphic
to $N_2\times D^2$. 
The diffeomorphism from $N_2\times D^2$ to $U$ can be chosen to be fiber
preserving: each fiber $q\times D^2$, $q\in N_2$, goes to a disk in some
plane $p\times\reals^2$, $p\in W$, where $p$ depends on $q$.
To be more precise, let us pick $U$ so that its intersection with every 
plane $p\times\reals^2$ is a disjoint union of a finite number
of two-dimensional disks. By means of the linear structure 
in $\reals^2$ we identify each of these disks with a small disk $D^2$ in 
$\reals^2$ centered at the origin. One can choose $U$ such that different 
fibers give rise to the same disk for all points of $N_2$. In this way, a 
neighborhood $U$ of $N_2$ in $W\times\reals^2$ turns into the direct 
product $N_2\times D^2$. 

Let $\tit$ and $\ty$ be the restrictions of $t$ and $y$ to the disk $D^2$ 
centered at the origin. We extend
$\tit$ and $\ty$ to $U$ using the direct product decomposition
$U=N_2\times D^2$ and denote by $t$ and $y$ the restrictions of
$t$ and $y$, as functions on $W\times\reals^2$, to $U$. Then 
$y=\ty$ and $t=\tit+h$, where $h$ is a function
on $N_2$. The restriction of $\sigma_W+dt\wedge dy$ to $U$ 
is therefore equal to $\sigma+d\tit\wedge d\ty+dh\wedge d\ty$. 

Denote by $\xi_h$ the Hamiltonian vector field of $h$ on $N_2$, i.e.,
$i_{\xi_h}\sigma=-dh$ and set $v=-\ty\xi_h$ on $U=N_2\times D^2$ by
using the direct product structure. A straightforward calculation
shows that the local flow $\psi^\tau$ of $v$ on $U$ is such that
$$
(\psi^\tau)^*(\sigma+d\tit\wedge d\ty+\tau dh\wedge d\ty)
=\sigma+d\tit\wedge d\ty 
$$
whenever $\psi^\tau$ is defined. Since $v$ vanishes on the section
$\ty=0$ of $N_2\times D^2$, the time-one flow $\psi$
is defined as a mapping of a small neighborhood of $\bar{N}_1$
to $U$.
Thus choosing $\epsilon>0$ and $\delta>0$ sufficiently small,
we get $\psi\colon N_1\times I\times J\to U$ with
$$
\psi^*(\sigma_W+dt\wedge dy)=\psi^*(\sigma+d\tit\wedge d\ty+dh\wedge d\ty)
=\sigma+d\tit\wedge d\ty
$$
and such that $\psi=id$ on $N_1\times I\times 0$. Also,
$\psi$ preserves $\tit$ and $\ty$ because $L_v\tit=L_v\ty=0$.
This completes the proof of the lemma. \qed

Denote by $F\colon \reals^2\to I\times J$ an arbitrary symplectic 
immersion. Without lost of generality we may assume that $F$ sends the
origin to the origin.
(To construct $F$ we may, for example, first squeeze symplectically
$\reals^2$ into a narrow infinite strip and then roll it up onto
an annulus lying in $I\times J$.) Let us extend the functions $f$ and
$g$ on $N$ to some smooth compactly supported functions on $N_1$,
which we denote by $f$ and $g$ again. 
As in the proof of Lemma (B') on p. 336 in \cite{gr:book},
define the embedding
$\tilde{j}\colon N_1\to N_1\times I\times J$ as the graph of 
$F\circ(f,g)\colon N_1\to I\times J$. Clearly, $\tilde{j}$ has compact support
and $\tilde{j}^*(\sigma+d\tit\wedge d\ty)=\sigma+df\wedge dg=\omega_{(k)}$.
The symplectic embedding 
$\tilde{i}\colon N\to W\times\reals^2$ such that 
$\tilde{i}^*(\sigma_W+dt\wedge dy)=\omega_{(k)}$ is just $\psi\circ\tilde{j}$
restricted to $N$.

{\em The second step} is to modify $\tilde{i}$ to make it fit
into $W\times\reals$. We will alter $\tilde{j}$ so that to transform
it into an embedding $j$ with the image in $N_1\times I\times 0$
and the same pull-back form. Then it will remain to set $i=\psi\circ j$.

Denote by $\pi$ the
projection $N_1\times I\times J\to N_1\times I$ along $J$.
We claim that the restriction of $\pi$ to the image
$\tilde{j}(N_1)$ is an embedding. 
This is an immediate consequence of the fact that $\tilde{j}(N_1)$ is a graph 
of some smooth function $N_1\to I\times J$. 

Thus $\tilde{j}(N_1)$ is a graph of a smooth function 
$\pi(\tilde{j}(N_1))\to J$ with compact support. Let us extend this 
function to a
smooth compactly supported function $H\colon N_1\times I\to J$.
Clearly, $\tilde{j}(N_1)$ lies on the graph $\Gamma$ of $H$. 

Finally, we claim that $\Gamma$ with the restriction of 
$\sigma+d\tit\wedge d\ty$ is ``symplectomorphic'' to $(N_1\times I, \sigma)$.
The symplectomorphism is given by the straightening of the 
characteristics on $\Gamma$. This is possible because
the $\tit$-component of a characteristic is nonzero and, on the
complement to a compact set, the characteristics are just the straight lines
parallel to $I$.

To give a more rigorous definition of the diffeomorphism,
let us think of $H$ as a time-dependent Hamiltonian on $N_1$ with
$I$ representing the time axis.
Denote its local flow in the extended phase-space $N_1\times I$ 
by $\phi^{\tit}$. (To be more precise about the notation, $\phi^{\tit}$ 
sends a point 
$(p,\tit')\in N_1\times I$ to the point $(\phi^{\tit}(p,\tit'),\tit+\tit')$.)
The characteristics of $\Gamma$ project under $\pi$ 
to the integral curves of $\phi^{\tit}$. (See, e.g., \cite{ar:math-meth}.)
Since $\supp H$ is compact in $N_1\times I$, there
exists $\tit_0>-\epsilon$ such that $\supp H$ is entirely contained in the 
region $\{\tit>\tit_0\}$. Now we define the desired symplectomorphism 
$\phi\colon \Gamma\to N_1\times I$ by sending
a point $(p, \tit, H(p,\tit))\in\Gamma$ to the point 
$(\phi^{\tit_0-\tit}(p, \tit), \tit)\in N_1\times I$ when $\tit>\tit_0$ 
and setting $\phi=id$ in the region $\tit_0\geq \tit>-\epsilon$.
By restricting the composition $\phi\circ\tilde{j}$ to $N$, we
obtain an embedding $j\colon N\to N_1\times I$ such that
$j^*\sigma=\tilde{j}^*(\sigma+d\tit\wedge d\ty)=\omega_{(k)}$.
Finally, $i_k=\psi\circ j$ satisfies the desired conditions: 
$i_k(N)\subset W\times\reals$ and $i_k^*\sigma_W=\omega_{(k)}$.

When $\p N=\emptyset$, the above argument goes through
for $N_2=N_1=N$ because the flow of $v$ exists for all times 
$\tau\in\reals$. Alternatively, in this case Theorem \ref{thm:embed}
follows immediately from Moser's theorem as is pointed out in Remark
\ref{rmk:moser1}.

\hfill\end{pf}

Before applying Theorem \ref{thm:embed}, let us recall the definition
of the horocycle flow. 

Let $\Sigma$ be a closed surface with a hyperbolic metric, i.e., a
metric with constant negative curvature $K=-1$. Denote by
$M=ST^*\Sigma$ the unit cotangent bundle of $\Sigma$ and by
$N$ a small closed tubular neighborhood of $M$ in $T^*\Sigma$. Let 
$\lambda$ be the restriction to $M$ (or $N$) of the canonical 
Liouville one-form ``$p\,dq$'' and by $\theta$ the connection form 
on the principle circle bundle $M\to \Sigma $ (taken with the negative sign). 
We also denote by $\theta$ its pull-back to $N$ under the radial 
projection $N\to M$. It is easy to see that $d\theta$ is the pull-back
of the metric area form on $\Sigma$.
For any $t\in [0,1]$, the sum $\omega_t=d\lambda+ td\theta$ is a 
symplectic form on $N$. 

In what follows, the pair $(M,\omega_1|_M)$ is referred to as the ``horocycle
flow'', for its characteristics form the flow lines of the true
horocycle flow. (See \cite{gi:Cambr} and \cite{gi:MathZ} for details.)
Our goal is to show that $(M,\omega_1|_M)$  can be embedded
``symplectically'' to $(\reals^5,\sigma)$ where $\sigma$ is the pull-back
to $\reals^5$ of the standard symplectic form on $\reals^4$. (By a
symplectic embedding we mean the one that pulls back $\sigma$ to $\omega_1$
and which is nowhere tangent to the characteristics of $\sigma$ in $\reals^5$.)
This will follow from (but is actually equivalent to) the existence of a 
symplectic embedding of $(N,\omega_1)$ to $(\reals^5,\sigma)$. 

By Theorem
\ref{thm:embed}, it is sufficient to find a symplectic embedding of
$(N,\omega_0)$, i.e., for the standard symplectic form on $N$. To construct
such an embedding, let us start with a Lagrangian immersion 
$\Sigma\to\reals^4$ (see \cite{gr:icm} and \cite{gr:book}), which we
may assume to have only simple double points. 
Extend this immersion to a symplectic immersion of a small neighborhood 
$V$ of $\Sigma$ in $T^*\Sigma$. The Lagrangian immersion $\Sigma\to\reals^4$
can be lifted to an embedding $\Sigma\to\reals^5$ so that $\sigma$ 
pulls back to the zero form. As a consequence, the symplectic embedding
$V\to \reals^4$ can be lifted to a symplectic embedding 
$(V,d\lambda)\to (\reals^5,\sigma)$. Applying a dilation in $\reals^4$
and in $T^*\Sigma$, we can
make $V$ arbitrarily large. In particular, we can take $V$ such that
$N\subset V$ which results into an embedding
$(N,\omega_0=d\lambda)\to (\reals^5,\sigma)$. Thus we have proved

\begin{cor}
\labell{cor:embedding}
There exists an embedding $j\colon N\to\reals^5$ such that 
$j^*\sigma=\omega_1$. There exists an embedding
$M\to\reals^5$ which is nowhere tangent to the characteristics of 
$\sigma$ in $\reals^5$ and such that the pull-back of $\sigma$ is $\omega_1$.
\end{cor}

\begin{rmk}
The existence of a symplectic embedding of $(N, \omega_t)$ into $\reals^6$ 
is an immediate consequence of Gromov's theorem on symplectic embeddings 
(\cite{gr:book}, pp. 335--336).
Theorem \ref{thm:embed} refines some of the assertions of Gromov's theorem 
by reducing the dimension of the ambient space by one. 
\end{rmk}

\begin{rmk}
It is not very clear how to show that $(M,\omega_1|_M)$ admits no symplectic 
embeddings into $\reals^4$. Assume that it does.
Then since $M$ is a closed hypersurface in $\reals^4$, the complement
$\reals^4\setminus M$ has two connected components: bounded, $U_b$, and
unbounded, $U_u$. Let $V_b$ and $V_u$, respectively, be the bounded and
unbounded components of $T^*\Sigma\setminus M$ equipped with $\omega_1$.
By a symplectic neighborhood theorem, one can attach either $V_b$ to $U_u$ 
(and $V_u$ to $U_b$) or $V_b$ to $U_b$ (and $V_u$ to $U_u$) with
a symplectic result. In the former case, the resulting manifold 
$V_b\cup_M U_u$ would
be symplectomorphic to $\reals^4$ at infinity, and so homeomorphic 
and (even symplectomorphic) to the standard $\reals^4$ by a theorem of Gromov 
and McDuff, \cite{gr:paper} and \cite{mcduff}. As a consequence, the
symplectic form on it must be exact. This is impossible because
$\Sigma\subset V_b$ and $\int_\Sigma\omega_1>0$. However,
the second way of attachment is not so easy to rule out.
\end{rmk}

The proofs of Theorems \ref{thm:main1} and \ref{thm:main2} repeat 
word-for-word the proofs of their higher-dimensional counterparts
(Theorems 2.1 and 2.5) given in \cite{gi:seifert}. In fact, the proof
of Theorem \ref{thm:main2} becomes even simpler than the proof of
Theorem 2.5  because now $M$ has codimension
one in $N$ and therefore one does not have the $D^{2k}$-component
($z$-coordinates). The condition (8), which the embedding
$M\to\reals^5$ must satisfy, is an immediate consequence of the fact
that the embedding is nowhere tangent to the characteristics or
equivalently that we have a symplectic embedding $N\to\reals^5$.
Finally, the argument used in \cite{gi:seifert} to prove Theorem 2.1
literally applies to Theorem \ref{thm:main1}.

\section{A list of Hamiltonian flows without periodic orbits} 
\labell{sec:list}

In this section we list all examples  (known to the author)
of smooth Hamiltonian systems on symplectic manifolds $(V,\omega)$ with a
compact regular energy level having no periodic orbits. 
This list divides in two, depending on whether or not the cohomology class 
of the symplectic form (near the level) is exact.

\bigskip

\noindent {\bf Case 1:}{\em~ The form $\omega$ is not required to be exact.}
In this case, one can take the torus $V=\TT^{2n}$, $2n\geq 4$, with an
irrational translation-invariant symplectic structure $\omega$. Then 
choose a Hamiltonian $H$ on $V$ so that the level $\{ H=1\}$ is the union 
of two standard embedded tori $\TT^{2n-1}\subset \TT^{2n}$. Since $\omega$
is irrational, the characteristics of $\omega|_{\TT^{2n-1}}$ form an
quasiperiodic flow on $\TT^{2n-1}$. Such a flow obviously has no periodic 
orbits. In fact, for a suitable choice of $H$, none of the levels
$\{H=c\}$ with $c\in (0.5, 1.5)$ carries a periodic orbit.
This example is due to Zehnder \cite{zehnder}. As shown by
M. Herman, the flow in question exhibits remarkable stability properties
\cite{herman1}, \cite{herman2}.

\bigskip

\noindent {\bf Case 2:}{\em~ The form $\omega$ is exact near the energy level.}
In this case we should distinguish whether $\dim V=4$ or $\dim V \geq 6$.

When $\dim V=4$, the only known example is the horocycle 
flow described as a Hamiltonian system in Section 
\ref{sec:proofs}. (Here $V=T^*\Sigma$, where $\Sigma$ is a compact surface 
with a metric of constant negative curvature $-1$,  $\omega=\omega_1$ is 
the twisted symplectic form, and $H$ is the standard metric Hamiltonian.) 
It is an old result of Hedlund that the horocycle flow 
has no periodic orbits on $M=\{H=1\}$ \cite{hedlund}. The 
horocycle flow naturally arises as a Hamiltonian system for the motion 
of a charge in a magnetic field on the surface $\Sigma$. This is the only 
known such system without periodic orbits. (See \cite{gi:Cambr} and 
\cite{gi:MathZ} for a detailed discussion.)

Observe that $N$ is
$G\backslash \mbox{\rm PSL}(2,\reals)\times (1-\epsilon,1+\epsilon)$ with
$G=\pi_1(\Sigma)$. Then $H$ becomes the projection to
the second component. The flow we just described is the Hamiltonian flow of
$H$ with respect to some symplectic form. Instead of $G=\pi_1(\Sigma)$ we
can take any discrete subgroup such that the quotient is compact and smooth.

An example of a $C^1$-smooth divergence-free vector field on $S^3$ having 
no periodic orbits is due G. Kuperberg, \cite{kug}, as well as the construction
of $C^\infty$-smooth volume preserving flows with a finite number of
periodic orbits on closed three-manifolds. These examples are
more along the lines of our Theorem \ref{thm:main2} rather than
genuine Hamiltonian flows.

Finally, when $\dim V\geq 6$, we have the flows on $\reals^{2n}$ obtained 
in the present paper and \cite{gi:seifert} by means of inserting a symplectic
plug to kill periodic orbits on a given energy level. When 
$2n\geq 8$, there is an alternative construction of such a
plug due to M. Herman \cite{herman-fax} as the symplectization
of Wilson's plug, \cite{wilson}. These methods allow one to use 
instead of an ellipsoid any compact hypersurface in 
$\reals^{2n}$ (with a finite number of closed characteristics) as the
initial manifold $Q$. For instance,
the methods apply to non simply connected hypersurfaces found by
Laudenbach \cite{laud}. Finally, one can combine these examples
by employing, when possible, one of the above manifolds as the core of the
plug to find new flows without closed orbits or with a finite number of them. 
(See Remark 3.6 (ii) of \cite{gi:seifert}.) For instance, taking $S^1$ as the 
core, Cieliebak \cite{ciel} constructed embeddings 
$S^{2n-1}\subset\reals^{2n}$, 
$2n\geq 4$, with a very interesting location and geometry of closed 
characteristics.


\begin{thebibliography}{ABC}

\bibitem[Ar]{ar:math-meth}
Arnold, V.I., {\em Mathematical methods of classical mechanics},
Graduate Texts in Math., {\bf 60}, Springer Verlag, New York, 1978.

\bibitem[Ci]{ciel}
Cieliebak, K., Symplectic boundaries: Creating and destroying closed 
characteristics, Preprint, 1995, to appear in {\em GAFA}.

\bibitem[Gi1]{gi:seifert}
Ginzburg, V.L., An embedding $S^{2n-1}\to\reals^{2n}$, $2n-1\geq 7$,
whose Hamiltonian flow has no periodic trajectories, 
{\em IMRN}, (1995), no. 2, 83--98.

\bibitem[Gi2]{gi:Cambr}
Ginzburg, V.L., On closed trajectories of a charge in a magnetic
field. An application of symplectic geometry, in {\em Contact
and Symplectic Geometry}, C.B. Thomas (Editor), Publications of
the Newton Institute, Cambridge University Press, Cambridge, 1996,
p. 131--148.

\bibitem[Gi3]{gi:MathZ}
Ginzburg, V.L., On the existence and non-existence of closed 
trajectories for some Hamiltonian flows, {\em Math. Z.},
{\bf 223} (1996), 397--409.

\bibitem[Gr1]{gr:icm}
Gromov, M., A topological technique for the construction of
solution of differential equations and inequalities, {\em ICM}.
 1970, Nice, vol. 2, pp. 221--225.

\bibitem[Gr2]{gr:paper}
Gromov, M., Pseudo holomorphic curves in symplectic manifolds,
{\em Invent. Math.}, {\bf 82} (1985, 307--347.

\bibitem[Gr3]{gr:book}
Gromov, M., {\em Partial differential relations},
Springer Verlag, New York, 1986.


\bibitem[Gh]{ghys}
Ghys, E., Construction de champs de vecteurs sans orbite
p\'{e}riodique (d'apr\`{e}s Krystyna Kuperberg), 
{\em S\'{e}m. Bourbaki}, 1993--1994, no 785, Juin 1994.

\bibitem[He]{hedlund}
Hedlund, G.A., Fuchsian groups and transitive horocycles,
{\em Duke Math. J.}, {\bf 2} (1936), 530--542.

\bibitem[Her1]{herman1}
Herman, M.-R., Examples de flots hamiltoniens dont aucune
perturbations en topologie $C^\infty$ n'a d'orbites 
p\'{e}riodiques sur ouvert de surfaces d'\'{e}nergies,
{\em C.R. Acad. Sci. Paris S\'{e}r. I Math.}, {\bf 312} 
(1991), 989--994.

\bibitem[Her2]{herman2}
Herman, M.-R., Diff\'{e}rentiabilit\'{e} optimale et contre-exemples \'{a} la
fermeture en topologie $C^\infty$ des orbites r\'{e}currentes de flots 
hamiltoniens. {\em C. R. Acad. Sci. Paris S\'{e}r. I Math.}, 
{\bf 313} (1991), 49--51. 

\bibitem[Her3]{herman-fax}
Herman, M.-R., Fax to Eliashberg.

\bibitem[HZ1]{ho-ze:per-sol}
Hofer, H., Zehnder, E., Periodic solution on hypersurfaces and a
result by C. Viterbo, {\em Invent. Math.}, {\bf 90} (1987), 1--9.

\bibitem[HZ2]{ho-ze:book}
Hofer, H., Zehnder, E., {\em Symplectic invariants and Hamiltonian
dynamics}, Birkh\"{a}user, Advanced Texts; Basel-Boston-Berlin, 1994.

\bibitem[KuG]{kug}
Kuperberg, G., A volume-preserving counterexample to the Seifert 
conjecture, {\em Comment. Math. Helv.}, {\bf 71} (1996), 70--97.

\bibitem[KuK]{kuk}
Kuperberg, K., A smooth counterexample to the Seifert conjecture in
dimension three, {\em Annals of Math.}, {\bf 140} (1994), 723--732.

\bibitem[La]{laud}
Laudenbach, F., Trois constructions an topologie symplectique,
Preprint, 1996, to appear in 
{\em Ann. Inst. H. Poincar\'{e}}.

\bibitem[McD]{mcduff}
McDuff, D., The structure of rational and rulled symplectic manifolds,
{\em Journ. Amer. Math. Soc.}, {\bf 3} (1990), 679--712; Erratum:
{\em Journ. Amer. Math. Soc.}, {\bf 5} (1992), 987--988.

\bibitem[Mo]{moser}
Moser, J., On the volume elements on a manifold, {\em Trans. Amer. Math.
Soc.}, {\bf 120} (1965), 286--294.

\bibitem[St]{str}
Struwe, M., Existence of periodic solutions of Hamiltonian
systems on almost every energy surfaces, {\em Bol. Soc. Bras. Mat.},
{\bf 20} (1990), 49--58.

\bibitem[Vi]{vi:thm}
Viterbo, C., A proof of Weinstein's conjecture in $\reals ^{2n}$,
{\em Ann. Inst. H. Poincar\'{e}, Anal. Non Lin\'{e}aire}, {\bf 4} (1987),
337--356.

\bibitem[We]{we:conj}
Weinstein, A., On the hypothesis of Rabinowitz' periodic orbit theorem,
{\em J. Diff. Eq.}, {\bf 33} (1979), 353--358.

\bibitem[Wi]{wilson}
Wilson, F., On the minimal sets of non-singular vector fields,
{\em Ann. Math.}, {\bf 84} (1966), 529--536.

\bibitem[Ze]{zehnder}
Zehnder, E., Remarks on periodic solutions on hypersurfaces, in
{\em Periodic solutions of hamiltonian systems and related topics}
by Rabinowitz et al., Reidel Publishing Co. (1987), 267--279.

\end{thebibliography}
\end{document}